\documentclass[12pt]{amsart}
\usepackage{amssymb}
\usepackage{amsbsy}
\usepackage{amscd}

\usepackage[mathscr]{eucal}

\DeclareSymbolFont{rsfs}{U}{rsfs}{m}{n}
\DeclareSymbolFontAlphabet{\mathrsfs}{rsfs}

\usepackage{a4wide}
\usepackage{verbatim}
\usepackage{version}
\usepackage{color}
\definecolor{darkspringgreen}{rgb}{0.09, 0.45, 0.27}
\definecolor{deepjunglegreen}{rgb}{0.0, 0.29, 0.29}
\newenvironment{NB}{
\color{red}{\bf NB}. \footnotesize
}{}
\newenvironment{NB2}{
\color{blue}{\bf NB2}. \footnotesize
}{}
\excludeversion{NB}
\excludeversion{NB2}
\excludeversion{NB3}
\makeatletter

\usepackage{url}
\usepackage{ifmtarg}

\usepackage{xr-hyper}
\usepackage[
pdftex,
bookmarks=true,
colorlinks=true,
citecolor=deepjunglegreen,
filecolor=deepjunglegreen,
debug=true,
unicode,
pdfnewwindow=true]{hyperref}
\usepackage{cleveref}

\crefname{Theorem}{Theorem}{Theorems}
\crefformat{Theorem}{Theorem~#2#1#3}
\crefname{section}{\S}{\S\S}
\crefformat{section}{\S#2#1#3} 
\crefmultiformat{section}{\S\S#2#1#3}{, #2#1#3}{, #2#1#3}{, #2#1#3}
\crefname{Lemma}{Lemma}{Lemmas}
\crefformat{Lemma}{Lemma~#2#1#3}
\crefname{Proposition}{Proposition}{Propositions}
\crefformat{Proposition}{Proposition~#2#1#3}
\crefname{Corollary}{Corollary}{Corollaries}
\crefformat{Corollary}{Corollary~#2#1#3}
\crefname{Definition}{Definition}{Definitions}
\crefformat{Definition}{Definition~#2#1#3}
\crefname{Remark}{Remark}{Remarks}
\crefformat{Remark}{Remark~#2#1#3}
\crefname{Remarks}{Remark}{Remarks}
\crefformat{Remarks}{Remark~#2#1#3}
\crefname{Conjecture}{Conjecture}{Conjectures}
\crefformat{Conjecture}{Conjecture~#2#1#3}
\crefname{figure}{Figure}{Figure}
\crefformat{figure}{Figure~#2#1#3}
\crefname{Example}{Example}{Example}
\crefformat{Example}{Example~#2#1#3}

\crefname{appendix}{Appendix}{Appendices}
\crefformat{appendix}{\S#2#1#3}
\crefformat{enumi}{#2#1#3}

\crefname{equation}{}{}

\crefname{Theorem*}{Theorem}{Theorems}
\crefformat{Theorem*}{Theorem~#2#1#3}

\usepackage{enumitem}
\usepackage{youngtab}
\usepackage{mathtools}
\usepackage{extpfeil}
\usepackage{multirow}
\usepackage{blkarray}
\usepackage{dynkin-diagrams}

\renewcommand{\thesubsection}{\thesection(\@roman\c@subsection)}
\newcounter{number}
\makeatother
\newtheorem{Corollary}[equation]{Corollary}

\newtheorem{Theorem*}{Theorem}

\theoremstyle{definition}

\newtheorem{Example}[equation]{Example}

\theoremstyle{remark}
\newtheorem{Remark}[equation]{Remark}

\numberwithin{equation}{section}

\newcommand{\defeq}{\overset{\operatorname{\scriptstyle def.}}{=}}
\newcommand{\CC}{{\mathbb C}}
\newcommand{\ZZ}{{\mathbb Z}}

\newcommand{\SL}{\operatorname{\rm SL}}

\newcommand{\GL}{\operatorname{GL}}

\newcommand{\algsl}{\operatorname{\mathfrak{sl}}} 

\newcommand{\so}{\operatorname{\mathfrak{so}}}

\newcommand{\g}{{\mathfrak g}}
\newcommand{\h}{{\mathfrak h}}

\newcommand{\Hom}{\operatorname{Hom}}

\renewcommand{\MR}[1]{}
\newcommand{\Wedge}{{\textstyle \bigwedge}}

\newcommand{\LL}{\mathbb L}

\newcommand{\tslabar}{\mathbin{
\setbox0=\hbox{/\!\!/\!\!/}\rule[0.4\ht0]{\wd0}{.3\dp0}\kern-\wd0\box0}}

\newcommand{\Gr}{\mathrm{Gr}}

\makeatletter
\newcommand{\cA}[1][{}]{%
  \@ifmtarg{#1}%
  {\mathcal A}
  {\mathcal A(#1)}
}
\newcommand{\cAh}[1][{}]{%
  \@ifmtarg{#1}%
  {\mathcal A_\hbar}
  {\mathcal A_\hbar(#1)}
}
\makeatother

\newcommand{\po}{\ar@{}[dr]|{\text{\pigpenfont R}}}
\newcommand{\pb}{\ar@{}[dr]|{\text{\pigpenfont J}}}
\newcommand{\pp}{\ar@{}[dr]|{\text{\pigpenfont P}}}

\newcommand{\cM}{\mathcal M}

\newcommand{\fg}{{\mathfrak{g}}}

\newcommand{\fM}{{\mathfrak{M}}}

\newcommand{\bv}{{\mathbf v}} 
\newcommand{\bw}{{\mathbf w}} 

\makeatletter
\newdimen\y@inside
\y@inside=\y@boxdim
\advance\y@inside by -.2\y@linethick
\definecolor{halfgray}{gray}{0.5}
\newcommand{\rf}{\color{halfgray}
\rule[-.2\y@inside]{\y@inside}{\y@inside}
}
\makeatother

\DeclareSymbolFont{symbolsC}{U}{pxsyc}{m}{n}
\SetSymbolFont{symbolsC}{bold}{U}{pxsyc}{bx}{n}
\DeclareFontSubstitution{U}{pxsyc}{m}{n}
\DeclareMathSymbol{\medcirc}{\mathbin}{symbolsC}{7}

\newcommand{\ch}{\operatorname{ch}}

\newcommand{\fin}{\mathrm{fin}}
\newcommand{\aff}{\mathrm{aff}}

\begin{document}

\title[Euler numbers of Hilbert schemes and quantum dimensions]
{Euler numbers of Hilbert schemes of points on
  simple surface singularities and quantum dimensions of standard
  modules of quantum affine algebras}

\author[H.~Nakajima]{Hiraku Nakajima}
\address{Kavli Institute for the Physics and Mathematics of the Universe (WPI),
  The University of Tokyo,
  5-1-5 Kashiwanoha, Kashiwa, Chiba, 277-8583,
  Japan
}
\email{hiraku.nakajima@ipmu.jp}
\address{Research Institute for Mathematical Sciences,
Kyoto University, Kyoto 606-8502,
Japan}

\begin{abstract}
  We prove the conjecture by Gyenge, N\'emethi and Szendr\H{o}i
  \cite{MR3671513,MR3799153} giving a formula of the generating
  function of Euler numbers of Hilbert schemes of points
  $\operatorname{Hilb}^n(\CC^2/\Gamma)$ on a simple singularity
  $\CC^2/\Gamma$, where $\Gamma$ is a finite subgroup of $\SL(2)$. We
  deduce it from the claim that quantum dimensions of standard modules
  for the quantum affine algebra associated with $\Gamma$ at
  $\zeta = \exp(\frac{2\pi i}{2(h^\vee+1)})$ are always $1$, which is
  a special case of a conjecture by Kuniba \cite{MR1202213}. Here
  $h^\vee$ is the dual Coxeter number. We also prove the claim, which
  was not known for $E_7$, $E_8$ before.
\end{abstract}

\maketitle

\section*{Introduction}

In this paper, we prove the conjecture by Gyenge, N\'emethi and
Szendr\H{o}i \cite{MR3671513,MR3799153} giving a formula of the
generating function of Euler numbers of Hilbert schemes of points
$\operatorname{Hilb}^n(\CC^2/\Gamma)$ on a simple singularity
$\CC^2/\Gamma$, where $\Gamma$ is a finite subgroup of $\SL(2)$. When
$\Gamma$ is of type $A$, Euler numbers were computed by Dijkgraaf and
Su\l kowski \cite{MR2391107}, and Toda \cite{MR3394393}. The formula in
\cite{MR3671513,MR3799153} is given in a different form, and makes
sense for arbitrary $\Gamma$. 
The formula was proved for type $D$, as well as type $A$, in
\cite{MR3799153}.
Our proof of the formula in type $E$ is new.

The formula \cite{MR3671513,MR3799153} is written as a specialization
of the character for the `Fock space', which is the tensor product of
the basic representation of the affine Lie algebra $\fg_{\aff}$
corresponding to $\Gamma$, and the usual Fock space representation of
the Heisenberg algebra. The Fock space was realized in a geometric way
as the direct sum of homology groups of various connected components
of $\operatorname{Hilb}^n(\CC^2)^\Gamma$ for all $n$. See
\cite{MR1989196}. The character is well-known (see e.g.,
\cite[\S12]{Kac})
\begin{equation*}
    Z =
    \prod_{m=1}^\infty
    (1 - e^{-m\delta})^{-(n+1)}
    \sum_{\vec{m}\in\ZZ^n} 
    e^{-\frac{(\vec{m}, \mathbf C\vec{m})}2\delta}
    \prod_{i=1}^n e^{-m_i\alpha_i},
\end{equation*}
where $\alpha_i$ is the $i$th simple root of $\fg_{\fin}$, the
underlying finite dimensional complex simple Lie algebra, $\delta$ is
the positive primitive imaginary root, $n$ is the rank of
$\fg_{\fin}$, and $\mathbf C$ is the finite Cartan matrix for
$\fg_{\fin}$. Here we omit the factor $e^{\Lambda_0}$ for the $0$th
fundamental weight $\Lambda_0$. The well-known formula for the basic
representation gives $(1 - e^{-m\delta})^{-n}$ instead of the power
$-(n+1)$. The additional $-1$ comes from the Heisenberg algebra.
\begin{NB}
  Suppose the dimension vector $\bv$ for $\fM(\bv,\Lambda_0)$ is
  $m\delta + \bv_{\fin}$ with
  $\bv_{\fin} = \sum_{i\neq 0} m_i \alpha_i$. By the Frenkel-Kac
  construction, we have the `bottom' case with
  $\dim \fM(\bv,\Lambda_0) = 0$, and raise its vacuum by the
  Heisenberg algebra. For the bottom, we have
  \begin{equation*}
    0 = \dim \fM(m\delta+\bv_{\fin},\Lambda_0)
    = 2m - (m\delta+\bv_{\fin}, \mathbf C(m\delta+\bv_{\fin}))
    = 2m - (\bv_{\fin}, \mathbf C \bv_{\fin})
  \end{equation*}
  as $\mathbf C\delta = 0$. Thus
  $m = \frac12 (\bv_{\fin}, \mathbf C \bv_{\fin})$.
\end{NB}%

The conjectural formula in \cite{MR3671513}, which we will prove, is
the following:
\begin{Theorem*}\label{thm:1}
  The generating function of Euler numbers of
  $\operatorname{Hilb}^n(\CC^2/\Gamma)$
  \begin{equation*}
    \sum_{n=0}^\infty e^{-n\delta} \chi(\operatorname{Hilb}^n(\CC^2/\Gamma))
  \end{equation*}
  is obtained from $Z$ by the substitution
  \begin{equation*}
    e^{-\alpha_i} = \exp\left( \frac{2\pi i}{h^\vee+1}\right)
    \quad (i\neq 0),
  \end{equation*}
  where $h^\vee$ is the dual Coxeter number.
\end{Theorem*}

Let us sketch the strategy of our proof.
In order to compute Euler numbers of
$\operatorname{Hilb}^n(\CC^2/\Gamma)$, we use a recent result of Craw,
Gammelgaard, Gyenge and Szendr\H{o}i \cite{2019arXiv191013420C}:
$\operatorname{Hilb}^n(\CC^2/\Gamma)$ with the reduced scheme
structure is isomorphic to a quiver variety
$\fM_{\zeta^\bullet}(n\delta,\Lambda_0)$, where $\zeta^\bullet$ is a
stability parameter
\begin{equation*}
  \zeta_i^\bullet = 0 \quad\text{for all $i\neq 0$},\quad
  \zeta_0^\bullet < 0.
\end{equation*}
This stability parameter is different from the standard one $\zeta$
\begin{equation*}
    \zeta_i < 0\quad\text{for all $i\in I$},
\end{equation*}
which gives connected components $\fM_{\zeta}(\bv, \Lambda_0)$ of
$\operatorname{Hilb}^{n}(\CC^2)^\Gamma$. Here $\bv$ is an isomorphism
class of $\CC[x,y]/I$ as a $\Gamma$-module for an ideal
$I \subset \CC[x,y]$ of the coordinate ring $\CC[x,y]$ of $\CC^2$,
invariant under the $\Gamma$-action.

As $\zeta^\bullet$ lies in walls given by roots in the space of
stability parameters, $\fM_{\zeta^\bullet}(n\delta,\Lambda_0)$ is
singular in general. Nevertheless it has a representation theoretic
meaning, as shown in \cite{Na-branching}: the quiver variety
$\fM_{\zeta^\bullet}(n\delta,\Lambda_0)$ is responsible for the
restriction of $\fg_{\aff}$ to $\fg_{\fin}$.
The essential geometric ingredient for this relation is a projective
morphism
\begin{equation*}
    \pi_{\zeta^\bullet,\zeta}\colon
    \fM_{\zeta}(\bv,\Lambda_0)
    \to \fM_{\zeta^\bullet}(n\delta,\Lambda_0)
    = \operatorname{Hilb}^n(\CC^2/\Gamma),
    \qquad n = \bv_0.
\end{equation*}
By a local description of singularities of quiver varieties
(\cite[\S2.7]{Na-branching}, which went back to
\cite[\S6]{Na-quiver}), fibers of $\pi_{\zeta^\bullet,\zeta}$ are
lagrangian subvarieties in quiver varieties $\fM_{\zeta}(\bv^s,\bw^s)$
associated with the finite $ADE$ quiver for $\fg_{\fin}$, the finite
dimensional complex simple Lie algebra underlying $\fg_{\aff}$.
As a simple application of this result, we observe that
\begin{itemize}
\item Euler numbers of $\operatorname{Hilb}^n(\CC^2/\Gamma)$ can be
  written in terms of Euler numbers of $\fM_\zeta(\bv,\Lambda_0)$ and
  quiver varieties $\fM_\zeta(\bv^s,\bw^s)$ of the finite $ADE$ type.
\end{itemize}
In fact, we get \emph{more} data than we need: we have a
stratification of $\operatorname{Hilb}^n(\CC^2/\Gamma)$ so that
$\pi_{\zeta^\bullet,\zeta}$ is a fiber bundle over each stratum. We
can compute Euler numbers of all strata from Euler numbers of
$\fM_\zeta(\bv,\Lambda_0)$ and quiver varieties
$\fM_\zeta(\bv^s,\bw^s)$ of the finite $ADE$ type. See \eqref{eq:5}
for the precise formula.

There are several algorithms\footnote{We make a distinction between an
  algorithm for a computation and an explicit computation as in
  \cite{Na:E8}. Namely a computation results finitely many $\pm$,
  $\times$, integers and variables. We do not require that the final
  expression is readable by a human. The expression like
  $\sum_{i=1}^{2^{(2^{100})}} a_i$ with explicit $a_i$ is an
  algorithm.} to compute Euler numbers of quiver varieties. For
example, \cite{MR2144973} and \cite{MR2651380}. They are complicated,
and hard to use in practice.
On the other hand, the conjectural formula in \cite{MR3671513} above
does \emph{not} contain such a complicated algorithm. It is just given
by a simple substitution.
It means that we should have a drastic simplification if we take a
\emph{linear combination} of complicated Euler numbers of quiver
varieties of finite $ADE$ type. We do not need to compute Euler
numbers of individual strata, as we only need their sum.

This simplification has a representation theoretic origin. Let us
explain it.

The above specialized character is called the \emph{quantum dimension}:
\begin{equation*}
  \dim_q V = \left.\ch V\right|_{e^{-\alpha_i}=\exp(\frac{2\pi i}{h^\vee+1})}.
\end{equation*}
Here $V$ is a finite dimensional representation $V$ of $\fg_{\fin}$.
It was introduced by Andersen~\cite[Def.~3.1]{MR1182414} (see also
Parshall-Wang~\cite{MR1247547}). We choose a specific root of unity
$\exp(\frac{2\pi i}{2(h^\vee+1)})$, where these papers study more
general roots of unity $\zeta$.

Therefore the conjectural formula in \cite{MR3671513} states that the
generating function of Euler numbers of
$\operatorname{Hilb}^n(\CC^2/\Gamma)$ is given by the quantum
dimension of the Fock space at $\exp(\frac{2\pi i}{2(h^\vee+1)})$,
restricted from $\fg_{\aff}$ to $\fg_{\fin}$, just keeping track
$e^{-\delta}$.

A direct sum of homology groups (which are isomorphic to complexified
$K$-group) of Lagrangian subvarieties of quiver varieties of the
finite $ADE$ type carries a structure of a finite dimensional
representation of the quantum loop algebra
$\mathbf U_{q}(\mathbf L\fg_{\fin})$. See \cite{Na-qaff}. It is called a
\emph{standard module} of $\mathbf U_{q}(\mathbf L\fg_{\fin})$.

It turns out that the above simplification is a consequence of the
following representation theoretic result of an independent interest:
\begin{Theorem*}\label{thm:2}
  The quantum dimention of arbitrary standard module of
  $\mathbf U_q(\mathbf L\fg_{\fin})$ of type $ADE$ at
  $\zeta=\exp(\frac{2\pi i}{2(h^\vee+1)})$ is equal to $1$.
\end{Theorem*}

Since standard modules are tensor products of $l$-fundamental modules
\cite{VV-std}, it is enough to prove this result for $l$-fundamental
modules.
The author is told by Naoi that this result is a special case of more
general conjecture posed by Kuniba \cite[Conj.~2 (A.6a)]{MR1202213}
(see also Kuniba-Nakanishi-Suzuki \cite[Conjecture~14.2]{MR2773889})
formulated for Kirillov-Reshetikhin modules.
(Recall $l$-fundamental modules are the simplest examples of
Kirillov-Reshetikhin modules.)
It is not difficult to check the general conjecture for type $A$. Type
$D$ case was shown in \cite{MR3043892}, while type $E_6$ case was
shown in \cite{MR3282650}.

Although we only need the \emph{simplest special} case of more general
conjecture, we could not find a proof of the relevant result for type
$E_7$, $E_8$ in the literature. Therefore we give its proof.
Fortunately the necessary explicit computation for type $E$ from the
algorithm in \cite{MR2144973} was already done in \cite{Na:E8} by
using a supercomputer. Alternatively we could quote the computation by
Kleber \cite{MR1436775}, which assumed fermionic formula conjectured
at that time. The fermionic formula was proved later by
Di~Francesco and Kedem~\cite{MR2428305}.
\begin{NB}
  based on many earlier works including \cite{MR1993360}, as well as
  results from cluster algebras.
\end{NB}%

We also give a proof of the known cases $A$, $D$, $E_6$ for the
completeness. We encounter a new feature in $E_7$, $E_8$, which did
not arise in other cases. Hence our check in the simplest case is yet
nontrivial.

Because it is about Kirillov-Reshetikhin modules, it suggests that a
general conjecture should be studied in the framework of cluster
algebras, as in \cite{MR2428305}. Note also that Euler numbers of
(graded) quiver varieties are understood in the context of cluster
algebras in a recent work of Bittmann \cite{2019arXiv191113110B}. Thus
the suggestion is compatible with the approach in this paper, though
cluster algebras play no role in this paper.

The paper is organized as follows. In \cref{sec:Euler} we deduce
\cref{thm:1} from \cref{thm:2} after recalling results from
\cite{Na-branching}. In \cref{sec:2} we prove \cref{thm:2}. It is
proved by the case by case analysis. In \cref{sec:rationally} we
discuss an additional topic, the rationally smoothness of
$\operatorname{Hilb}^n(\CC^2/\Gamma)$. It is rationally smooth if
$n=2$, but not so in general.

\subsection*{Acknowledgments}

The author thanks
A.~Gyenge
and
B.~Szendr\H{o}i
for discussion on \cite{MR3671513,MR3799153,2019arXiv191013420C}, and
A.~Craw
on \cite{2019arXiv191013420C} and \cite{MR1202213,MR2773889}.
He also thanks
A.~Kuniba,
S.~Naito,
T.~Nakanishi
and
K.~Naoi
for discussion on the conjecture in \cite{MR1202213,MR2773889},
D.~Muthiah
for discussion on rational smoothness,
and
R.~Yamagishi
for discussion on \cite{2017arXiv170905886Y}.

The research is supported in part by the World Premier International
Research Center Initiative (WPI Initiative), MEXT, Japan, and by JSPS
Grants Numbers 16H06335, 19K21828.

\section{Euler numbers of quiver varieties}\label{sec:Euler}

\subsection{Quiver varieties}

Let $\Gamma$ be a nontrivial finite subgroup of $\SL(2)$.
We define the affine Dynkin diagram via the McKay correspondence (see
e.g., \cite[Ch.~8]{MR3526103} for detail): let
$\{ \rho_i\}_{i\in I}$ be the set of isomorphism classes of
irreducible representations of $\Gamma$ with the trivial
representation $\rho_0$. We identify $\rho_i$ with a vertex of a
graph. We draw $a_{ij}$ edges between $\rho_i$ and $\rho_j$ where
$a_{ij} = \dim \Hom_\Gamma(\rho_i, \rho\otimes\rho_j) =
\dim\Hom_\Gamma(\rho_j,\rho\otimes\rho_i)$, where $\rho$ is the
$2$-dimensional representation of $\Gamma$ given by the inclusion
$\Gamma\subset\SL(2)$. Then the graph is an affine Dynkin diagram of
type $ADE$. Let $\fg_{\aff}$ denote the corresponding affine Lie
algebra, and $\fg_{\fin}$ the underlying finite dimensional complex
simple Lie algebra corresponding to the Dynkin diagram obtained from
the affine one by removing $\rho_0$. Let $n$ be the rank of
$\fg_{\fin}$, which is the number of vertices in $I$ minus $1$.

We use the convention of the root system for $\fg_{\aff}$ as in
\cite[Ch.~6 and \S12.4]{Kac}.
Let $\alpha_i$ be the $i$th simple root of $\fg_{\aff}$ corresponding
to $\rho_i$, $\delta$ be the primitive positive imaginary root of
$\fg_{\aff}$. We have $\delta = \sum a_i \alpha_i$, and $a_i$ is equal
to the dimension of $\rho_i$.
Let $\alpha_i^\vee$ be the $i$th simple coroot. We take the scaling
element $d$ satisfying
\begin{equation*}
  \langle\alpha_i, d\rangle = \delta_{i0}.
\end{equation*}
Then $\langle \alpha_i^\vee, d\rangle_{i\in I}$ is a base of
$\mathfrak h_{\aff}$, the Cartan subalgebra of $\fg_{\aff}$.
We define the fundamental weights $\Lambda_i$ ($i\in I$) by
\begin{equation*}
  \langle \Lambda_i, \alpha_j^\vee\rangle = \delta_{ij},\quad
  \langle \Lambda_i, d\rangle = 0.
\end{equation*}
Then $\langle \alpha_i, \Lambda_0\rangle_{i\in I}$ forms a base of
$\mathfrak h_{\aff}^*$.

We choose an orientation of edges in the affine Dynkin diagram and
consider the corresponding affine quiver $Q = (I,\Omega)$.
We take dimension vectors $\bw = (\bw_i)$, and
$\bv = (\bv_i)\in \ZZ_{\ge 0}^{I}$, and
consider quiver varieties $\fM_\zeta(\bv,\bw)$,
$\fM_{\zeta^\bullet}(\bv,\bw)$, where $\zeta$, $\zeta^\bullet$ are
stability parameters such that
\begin{equation*}
  \zeta_i < 0\quad\text{for all $i\in I$},\qquad
  \zeta_i^\bullet = 0 \quad\text{for all $i\neq 0$},\quad
  \zeta_0^\bullet < 0.
\end{equation*}
See \cite[\S2]{Na-branching} for the definition of quiver varieties
for these stability conditions.
Since $\zeta^\bullet$ lives in the boundary of a chamber containing
$\zeta$, we have a projective morphism
\begin{equation*}
  \pi_{\zeta^\bullet,\zeta}\colon
  \fM_\zeta(\bv,\bw)\to \fM_{\zeta^\bullet}(\bv,\bw).
\end{equation*}
See \cite[\S2]{Na-branching}.
In fact, this is an example studied in \cite[\S2.8]{Na-branching}
associated with a division $I = I^0\sqcup I^+$ with
$I^0 = I\setminus \{0\}$, $I^+ = \{0\}$.

We identify $\bv$, $\bw$ with weights of $\fg_{\aff}$ by
\begin{equation*}
  \bv = \sum \bv_i \alpha_i, \quad \bw = \sum \bw_i \Lambda_i.
\end{equation*}
It is convenient to use a different convention for the dimension
vector $\bv$:
\begin{equation}\label{eq:3}
  \bv = m\delta + \sum_{i\in I^0} m_i \alpha_i, \qquad \text{i.e., }
  m = \bv_0, \quad m_i = \bv_i - \bv_0 a_i.
\end{equation}

Let us take $\bw = \Lambda_0$.
It is well-known that $\fM_\zeta(\bv,\Lambda_0)$ is the $\Gamma$-fixed
point component of Hilbert schemes $I$ of points in the affine plane
$\CC^2$ such that $\CC[x,y]/I$ is isomorphic to
$\bigoplus \rho_i^{\oplus \bv_i}$, as a $\Gamma$-module.
Euler number for $\fM_\zeta(\bv,\Lambda_0)$ is also known. It was
given in \cite{MR1989196}. Since it was stated without an explanation,
let us explain how it is derived. We change the stability condition
$\zeta$ to see that $\fM_\zeta(\bv,\bw)$ is diffeomorphic to a moduli
space of framed rank $1$ torsion free sheaves on the minimal
resolution of $\CC^2/\Gamma$. Then rank $1$ torsion free sheaves are
ideal sheaves twisted by line bundles, hence Euler numbers are given
by G\"ottsche formula for Hilbert schemes of points \cite{MR1032930}.
Moreover this latter picture gives the Frenkel-Kac construction of the
basic representation of $\fg_{\aff}$ (see e.g., \cite[\S14.8]{Kac}),
hence
we get the formula of $Z$ in Introduction. We have used the convention
in \eqref{eq:3}.

We can define a structure of a representation of the affine Lie
algebra $\fg_{\aff}$ on the direct sum of homology groups
$\bigoplus_\bv H_*(\fM_\zeta(\bv,\Lambda_0))$. See \cite{MR1989196}
and references therein. We can also construct a structure of a
representation of the Heisenberg algebra commuting with $\fg_{\aff}$
by \cite{MR1441880,Lecture}.

On the other hand, Craw, Gammelgaard, Gyenge and Szendr\H{o}i
\cite{2019arXiv191013420C} recently proved that the Hilbert scheme
$\operatorname{Hilb}^n(\CC^2/\Gamma)$ of $n$ points in $\CC^2/\Gamma$
with the reduced scheme structure is isomorphic to
$\fM_{\zeta^\bullet}(\bv,\Lambda_0)$ with $\bv = n\delta$.

\subsection{Stratification}\label{subsec:strat}

As is explained in \cite[\S2]{Na-branching},
$\fM_{\zeta^\bullet}(\bv,\bw)$ parametrizes $S$-equivalence classes of
$\zeta^\bullet$-semistable framed representations of the preprojective
algebra. Therefore its points are represented by direct sum of
$\zeta^\bullet$-stable representations. Under the above choice of
$\zeta^\bullet$, we have one distinguished summand, giving a point in
$\fM_{\zeta^\bullet}^{\mathrm{s}}(\bv',\bw)$ for some $\bv'\le\bv$
(component-wise) such that $\bv'_0 = \bv_0$, and other summands are
simple representations $S_i$ with $i\in I^0$. See
\cite[\S2.6]{Na-branching}.
Since multiplicities of $S_i$ can be read off from the difference
$\bv - \bv'$, we can regard
$\fM_{\zeta^\bullet}^{\mathrm{s}}(\bv',\bw)$ as a subset of
$\fM_{\zeta^\bullet}(\bv,\bw)$.
By \cite[Prop.~2.30]{Na-branching}
$\fM^{\mathrm{s}}_{\zeta^\bullet}(\bv',\bw)\neq\emptyset$ if and only
if $\bw - \bv'$ is an $I^0$-dominant weight appearing in the basic
representation $V(\Lambda_0)$ of the affine Lie algebra $\fg$
associated with $\Gamma$.
Here $I^0$-dominant means that
$\langle \bw - \bv', \alpha_i^\vee\rangle \ge 0$ for $i\in I^0$.
We thus have the stratification
\begin{equation*}
  \fM_{\zeta^\bullet}(\bv,\bw)
  = \bigsqcup_{\text{$\bv'$ as above}} \fM_{\zeta^\bullet}^{\mathrm{s}}(\bv',\bw). 
\end{equation*}

Moreover the transversal slice to the stratum
$\fM^{\mathrm{s}}_{\zeta^\bullet}(\bv',\bw)$ is locally isomorphic to
a quiver variety $\fM_0(\bv^s,\bw^s)$ around $0$, associated with the
finite $ADE$ quiver $Q\setminus \{\rho_0\}$ such that dimension
vectors are given by
\begin{equation*}
  \bv^s = \bv - \bv', \qquad
  \bw^s_i = \langle \bw - \bv', \alpha_i^\vee\rangle \quad i\in I^0.
\end{equation*}
Note that $\bv-\bv'$ has no $0$th component as $\bv'_0 = \bv_0$.
Note also that $\bw^s_i \ge 0$, as $\bw - \bv'$ is $I^0$-dominant.
It is also known that the inverse image of the slice in
$\fM_\zeta(\bv,\bw)$ under $\pi_{\zeta^\bullet,\zeta}$ is locally
isomorphic to $\fM_{\zeta}(\bv^s,\bw^s)$ around
$\mathfrak L(\bv^s,\bw^s)$, the inverse image of the origin $0$ of
$\fM_0(\bv^s,\bw^s)$ under the projective morphism
$\fM_\zeta(\bv^s,\bw^s)\to \fM_0(\bv^s,\bw^s)$.
See \cite[\S2.7]{Na-branching} and the references therein for these
claims on transversal slices. There is also an algebraic approach in
\cite{CB:normal}.

For $\bw = \Lambda_0$, we have
\begin{equation*}
  \bw^s_i = - \sum_{j\in I^0} m'_j \langle \alpha_j, \alpha_i^\vee\rangle,
  \qquad \text{where $
    \bv' = m\delta + \sum_{j\in I^0} m'_j \alpha_j$
    in the convention \eqref{eq:3}}.
\end{equation*}
Note that $\bv = m\delta + \sum m_i \alpha_i$ and $\bv'$ share the
same $m$ for the coefficient of $\delta$, as $\bv_0 = \bv'_0$.
In particular,
\begin{equation}\label{eq:6}
    \sum_{i\in I^0} \bw^s_i \Lambda_i - \bv^s_i \alpha_i
    \begin{NB}
        = - \sum_{i\in I^0} \sum_{j\in I^0} m'_j \langle\alpha_j,
        \alpha_i^\vee \rangle \Lambda_i - (m_i - m'_i)\alpha_i
    \end{NB}
    = - \sum_{i\in I^0} m_i \alpha_i.
\end{equation}
We also note
\begin{equation}\label{eq:7}
    m'_i \le 0 \qquad (i\in I^0),
\end{equation}
as $(m'_i)_{i\in I^0} = - (\bw^s)_{i\in I^0} \mathbf C^{-1}$, and the
inverse of the Cartan matrix $\mathbf C$ has positive entries. This
was proved in \cite[Prop.~A.1]{2019arXiv191013420C} by a different
method.

\subsection{Examples}

\begin{Example}
  Let
  $\fM_{\zeta^\bullet}(n\delta, \Lambda_0)\cong
  \operatorname{Hilb}^n(\CC^2/\Gamma)$.
  Consider a stratum containing $(n-1)$-distinct points in
  $\CC^2\setminus\{0\}/\Gamma$ together with trivial
  representations. We have
  $\bv' = (n-1)\delta + \alpha_0 = n\delta - \sum_{i\in I^0} a_i
  \alpha_i$, and hence
  \begin{equation*}
    \bv^s = \sum_{i\in I^0} a_i \alpha_i, \quad
    \bw^s_i = \sum_{j\in I^0} C_{ij} a_j.
\end{equation*}
Note that $\bw^s$ has entries $1$ at vertices in the finite quiver
$I^0$ which are connected to the $0$-vertex in the affine quiver. In
particular, $\fM_0(\bv^s,\bw^s)$ is $\CC^2/\Gamma$,
$\fM_{\zeta}(\bv^s,\bw^s)$ is its minimal resolution, the very
first example of a quiver variety considered by Kronheimer \cite{Kr},
before the definition of quiver varieties was introduced.

This is obvious, as the transversal slice is
$\fM_{\zeta^\bullet}(\delta,\Lambda_0) \cong \CC^2/\Gamma$, as we can
ignore $(n-1)$-distinct points.
\end{Example}

\begin{Example}[Yamagishi \cite{2017arXiv170905886Y}]\label{ex:Yam}
  Consider
  $\fM_{\zeta^\bullet}(2\delta, \Lambda_0) =
  \operatorname{Hilb}^2(\CC^2/\Gamma)$. The formal neighborhoods of
  fibers of $\pi_{\zeta^\bullet,\zeta}$ over $0$-dimensional strata in
  $\fM_{\zeta^\bullet}(2\delta,\Lambda_0)$ were determined by Yamagishi
  \cite{2017arXiv170905886Y}. He identified the formal neighborhood
  with that of the intersection of the nilpotent cone for the complex
  simple Lie algebra $\fg_{\fin}$ and the Slodowy slice to a
  `sub-subregular' orbit. There is only one $0$-dimensional stratum
  except type $A_1$, $A_2$, $D_n$, while there are none for $A_1$,
  $A_2$, two for $D_n$ ($n > 4$) and three for $D_4$.
  We can give corresponding vector $\bv'$ as follows:
  \begin{itemize}
  \item If $Q$ is of type $A_1^{(1)}$ or $A_2^{(1)}$, there is no such
    $\bv'$.
  \item If $Q$ is of type $A_{n}^{(1)}$ with $n > 2$,
    $\bv' = 2\alpha_0 + \alpha_1 + \alpha_{n}$.
  \item If $Q$ is of type $D_n^{(1)}$ with $n > 4$,
    $\bv' = \dynkin[edge
    length=.6cm,labels={2,0,2,2,2,2,1,1}]{D}[1]{***.****} = 2\alpha_0
    + 2\alpha_2 + \dots + 2\alpha_{n-2} + \alpha_{n-1} + \alpha_n$ or
    $\bv' = \dynkin[edge
    length=.6cm,labels={2,1,2,1,0,0,0,0,0}]{D}[1]{****.****} =
    2\alpha_0 + \alpha_1 + 2\alpha_2 + \alpha_3$.
  \item If $Q$ is of type $D_4^{(1)}$, we have three possibilities: in
    addition to the second example in $D_n^{(1)}$ above, we have
    $(\alpha_1,\alpha_3) \mapsto (\alpha_1,\alpha_4)$,
    $(\alpha_3,\alpha_4)$, as its cyclic permutation.

  \item If $Q$ is of type $E_8^{(1)}$, we have
    $\bv' = \dynkin[edge
    length=.6cm,labels={2,0,1,1,2,2,2,2,2}]{E}[1]{8}$. The cases
    $E_6^{(1)}$, $E_7^{(1)}$ are similar.
  \end{itemize}
  Here $0$th vertex is $\circ$, and other vertices are $\bullet$.

  Transversal slices are
  \begin{itemize}
  \item If $Q$ is of type $A_{n}^{(1)}$ with $n>2$,
    $\bv^s = \dynkin[edge
    length=.6cm,labels={1,2,2,2,2,1}]{A}{***.***}$,
    $\bw^s = \dynkin[edge
    length=.6cm,labels={0,1,0,0,1,0}]{A}{***.***}$.

  \item If $Q$ is of type $D_n^{(1)}$,
    $\bv^s = \dynkin[edge
    length=.6cm,labels={2,2,2,2,1,1}]{D}{***.***}$,
    $\bw^s = \dynkin[edge
    length=.6cm,labels={2,0,0,0,0,0}]{D}{***.***}$ in the first case,
    and
    $\bv^s = \dynkin[edge
    length=.6cm,labels={1,2,3,4,4,4,4,2,2}]{D}{*****.****}$,
    $\bw^s = \dynkin[edge
    length=.6cm,labels={0,0,0,1,0,0,0,0,0}]{D}{*****.****}$ in the
    second case.

  \item If $Q$ is of type $D_4^{(1)}$,
    $\bv^s = \dynkin[edge length=.6cm,labels={1,2,1,2}]{D}{4}$,
    $\bw^s = \dynkin[edge length=.6cm,labels={0,0,0,2}]{D}{4}$, and
    its cyclic permutation.

  \item If $Q$ is of type $E_8^{(1)}$,
    $\bv^s = \dynkin[edge
    length=.6cm,labels={4,5,7,10,8,6,4,2}]{E}{8}$,
    $\bw^s = \dynkin[edge
    length=.6cm,labels={1,0,0,0,0,0,0,0}]{E}{8}$.
\end{itemize}

In particular, it means that the quiver variety $\fM_0(\bv^s,\bw^s)$
is isomorphic to the intersection of Slodowy slice and the nilpotent
cone in the formal neighborhood of the origin.

This result was known before for type $A_n$ and the second case of
$D_n$. The type $A_n$ case was proved in \cite[\S8]{Na-quiver}. See
also \cite{MR2130242}. For the second case of type $D_n$, it was shown
in \cite{MR3269179}. In these cases, $\fM_0(\bv^s,\bw^s)$ itself is
isomorphic to the intersection of Slodowy slice and the nilpotent
cone, not only in the formal neighborhood.
The exceptional cases are conjectured, but not shown as far as the
author knows.
\end{Example}

\begin{Remark}
  In cases for the above Example, it is known that corresponding
  Coulomb branches of the quiver gauge theories (which are affine
  Grassmannian slices by \cite{2016arXiv160403625B}) are
  `next-to-minimal' nilpotent orbits closures in $\fg_{\fin}$ by
  \cite{MR3131493}.
\end{Remark}

\begin{Example}[Yamagishi \cite{2017arXiv170905886Y}]
  Let us consider $\bv = \delta+\alpha_0$, $\bw = \Lambda_0$, hence
  $\fM_{\zeta^\bullet}(\delta+\alpha_0,\Lambda_0)$. Note that this is
  different from
  $\fM_{\zeta^\bullet}(\delta,\Lambda_0) =
  \operatorname{Hilb}^1(\CC^2/\Gamma) = \CC^2/\Gamma$. Let us show
  that $\fM_{\zeta^\bullet}(\delta+\alpha_0,\Lambda_0)$ is \emph{not}
  isomorphic to $\CC^2/\Gamma$.

  This $\fM_{\zeta^\bullet}(\delta+\alpha_0,\Lambda_0)$ appears as the
  union of strata in
  $\fM_{\zeta^\bullet}(2\delta,\Lambda_0) =
  \operatorname{Hilb}^2(\CC^2/\Gamma)$. Thus all possible strata are
  determined already as above: besides the open stratum for
  $\bv' = \bv$, we have a single $0$-dimensional stratum except for
  $A_1$, $A_2$, $D_n$, none for $A_1$, $A_2$, three for $D_4$, two for
  $D_n$ ($n > 4$).
  The transversal slice $\fM_0(\bv^s,\bw^s)$ is given by the same
  $\bw^s$ as above, $\bv^s$ is subtracting
  $\sum_{i\in I^0} a_i\alpha_i$ from the above example. Concretely it
  is
  \begin{itemize}
  \item For type $A_n^{(1)}$ with $n > 2$,
    $\bv^s = \dynkin[edge
    length=.6cm,labels={0,1,1,1,1,0}]{A}{***.***}$,
    $\bw^s = \dynkin[edge
    length=.6cm,labels={0,1,0,0,1,0}]{A}{***.***}$.

  \item For type $D_n^{(1)}$ with $n > 4$,
    $\bv^s = \dynkin[edge
    length=.6cm,labels={1,0,0,0,0,0}]{D}{***.***}$,
    $\bw^s = \dynkin[edge
    length=.6cm,labels={2,0,0,0,0,0}]{D}{***.***}$ in the first case,
    and
    $\bv^s = \dynkin[edge
    length=.6cm,labels={0,0,1,2,2,2,2,1,1}]{D}{*****.****}$,
    $\bw^s = \dynkin[edge
    length=.6cm,labels={0,0,0,1,0,0,0,0,0}]{D}{*****.****}$ in the
    second case.

  \item For type $D_4^{(1)}$,
    $\bv^s = \dynkin[edge length=.6cm,labels={0,0,0,1}]{D}{4}$,
    $\bw^s = \dynkin[edge length=.6cm,labels={0,0,0,2}]{D}{4}$, and
    its cyclic permutation.

  \item For type $E_8^{(1)}$,
    $\bv^s = \dynkin[edge
    length=.6cm,labels={2,2,3,4,3,2,1,0}]{E}{8}$,
    $\bw^s = \dynkin[edge
    length=.6cm,labels={1,0,0,0,0,0,0,0}]{E}{8}$.
  \end{itemize}

  Let us consider
  \(
    \pi_{\zeta^\bullet,\zeta}\colon \fM_{\zeta}(\bv,\bw)
    \to \fM_{\zeta^\bullet}(\bv,\bw).
  \)
  The domain $\fM_{\zeta}(\bv,\bw)$ is the minimal resolution of
  $\CC^2/\Gamma$. The above description of the transversal slice
  immediately conclude the following:
  \begin{itemize}
  \item $\fM_{\zeta^\bullet}(\bv,\bw)$ is obtained from the minimal
    resolution of $\CC^2/\Gamma$ by collapsing $\mathbb P^1$'s whose
    vertices are not connected to the $0$-vertex in the affine Dynkin
    diagram.
  \end{itemize}

  There are no such vertices for $A_1$, $A_2$. For type $A_n$
  ($n > 2$), we collapse $(n-2)$ $\mathbb P^1$'s corresponding to
  vertices except the leftmost and rightmost. For type $D_n$, we
  collapse $(n-1)$ $\mathbb P^1$'s except the second one from the
  left, and producing singularities of type $A_1$ and $D_{n-2}$, where
  we understand $D_2 = A_1\times A_1$, $D_3 = A_3$.
  For type $E_8$, the $7$ $\mathbb P^1$'s except the rightmost one are
  collapsed.
\end{Example}

\subsection{Euler numbers}

By \cref{subsec:strat} we can relate Euler numbers of
$\fM_{\zeta}(\bv,\bw)$, $\fM_{\zeta^\bullet}(\bv',\bw)$ and
$\fM_{\zeta}(\bv^s,\bw^s)$. Let $\chi(\ )$ denote the Euler number of
a space. We have
\begin{equation*}
  \begin{split}
  \chi(\fM_{\zeta}(\bv,\bw))
  & = \sum_{\bv'} \chi(\fM_{\zeta^\bullet}^{\mathrm{s}}(\bv',\bw))
  \chi(\mathfrak L_{\zeta}(\bv^s,\bw^s)) \\
  & = \sum_{\bv'} \chi(\fM_{\zeta^\bullet}^{\mathrm{s}}(\bv',\bw))
  \chi(\fM_{\zeta}(\bv^s,\bw^s)).
  \end{split}
\end{equation*}
It is known \cite[Cor.~5.5]{Na-quiver} that the central fiber
$\mathfrak L_{\zeta}(\bv^s,\bw^s)$ of
$\fM_{\zeta}(\bv^s,\bw^s)\to \fM_0(\bv^s,\bw^s)$ is homotopic to
$\fM_{\zeta}(\bv^s,\bw^s)$, hence the second equality follows.

From $\chi(\fM_{\zeta}(\bv,\bw))$,
$\chi(\fM_{\zeta}(\bv^s,\bw^s))$, we compute
$\chi(\fM_{\zeta^\bullet}^{\mathrm{s}}(\bv',\bw))$ recursively as
\begin{equation}\label{eq:5}
  \chi(\fM_{\zeta^\bullet}^{\mathrm{s}}(\bv,\bw))
  = \chi(\fM_\zeta(\bv,\bw)) - \sum_{\bv'\neq \bv}
  \chi(\fM_{\zeta^\bullet}^{\mathrm{s}}(\bv',\bw))
  \chi(\fM_{\zeta}(\bv^s,\bw^s)).
\end{equation}
Here we use $\bv^s = 0$, hence $\fM_{\zeta}(\bv^s,\bw^s)$ is a point
for $\bv'=\bv$.

We take the generating function of Euler numbers as
\begin{equation}\label{eq:8}
  \sum_{\bv} \chi(\fM_{\zeta}(\bv,\bw)) e^{-\bv}
  = \sum_{\bv'} \chi(\fM_{\zeta^\bullet}^{\mathrm{s}}(\bv',\bw)) e^{-\bv'}
  \sum_{\bv^s}
  \chi(\fM_{\zeta}(\bv^s,\bw-\bv'|_{I^0})) e^{-\bv^s},
\end{equation}
where
$\bw - \bv'|_{I^0} = \bw^s = \sum_{i\in I^0} \langle \bw -
\bv',\alpha_i^\vee\rangle \Lambda_i$.

We claim
\begin{equation}\label{eq:1}
  \left.\sum_{\bv^s}
    \chi(\fM_{\zeta}(\bv^s,\bw^s)) 
    \prod_{i\in I^0} e^{\bw^s_i \Lambda_i
      - \bv^s_i \alpha_i}
  \right|_{
    e^{-\alpha_i} = \exp(\frac{2\pi i}{h^\vee+1})
  }
  = 1.
\end{equation}
\begin{NB}
  $\dim \fM_{\zeta}(\bv^s,\bw^s) = \langle 2\bw^s-\bv^s, \bv^s\rangle$.
\end{NB}%

We take $\bw = \Lambda_0$ and switch to the convention
\eqref{eq:3}. Then the left hand side of \eqref{eq:8} is $Z$ in
Introduction. By \eqref{eq:6}, \eqref{eq:1} implies that
\begin{NB}
\begin{multline*}
    \sum_{m, \vec{m}} \chi(\fM_{\zeta}(m\delta+\sum_{i\in I^0} m_i\alpha_i,
  \Lambda_0)) 
  e^{-m\delta} \prod_{i\in I^0} e^{-m_i\alpha_i}
\\
  = \sum_{m,\vec{m}'} \chi(\fM_{\zeta^\bullet}^{\mathrm{s}}(
  m\delta+\sum_{i\in I^0} m'_i\alpha_i,\Lambda_0)) 
  e^{-m\delta} \prod_{i\in I^0} e^{-m_i\alpha_i}
  \sum_{\vec{m''}}
  \chi(\fM_{\zeta}(\sum_{i\in I^0} m''_i\alpha_i,\bw^s)) 
  \prod_{i\in I^0}e^{-m''_i\alpha_i}.
\end{multline*}
\end{NB}
\begin{equation*}
    \left.Z\right|_{
      e^{-\alpha_i} = \exp(\frac{2\pi i}{h^\vee+1})
    }
    = \sum_{m,\vec{m}'} \chi(\fM_{\zeta^\bullet}^{\mathrm{s}}(
  m\delta+\sum_{i\in I^0} m_i'\alpha_i,\Lambda_0)) 
  e^{-m\delta}.
\end{equation*}
By \eqref{eq:7} we consider
$\fM_{\zeta^\bullet}^{\mathrm{s}}( m\delta+\sum_{i\in I^0}
m_i'\alpha_i,\Lambda_0)$ as a stratum of
$\fM_{\zeta^\bullet}(m\delta,\Lambda_0) =
\operatorname{Hilb}^m(\CC^2/\Gamma)$. Therefore the right hand side
of the above is the generating function of Euler numbers of
$\operatorname{Hilb}^m(\CC^2/\Gamma)$.
Thus we have proved \cref{thm:1}.

Recall that
$\bigoplus_{\bv^s} H_*(\mathfrak L_{\zeta}(\bv^s,\bw^s)) \cong
\bigoplus_{\bv^s} K(\mathfrak L_{\zeta}(\bv^s,\bw^s))\otimes_{\ZZ}\CC$
is the so-called standard module of the quantum loop algebra
$\mathbf U_q(\mathbf L\fg_{\fin})$, specialized at $q=1$
\cite{Na-qaff}.
The above \eqref{eq:1} means that its character (as a
$\fg_{\fin}$-module), specialized at
$e^{-\alpha_i} = \exp(\frac{2\pi i}{h^\vee+1})$ is equal to $1$.
As we mentioned in Introduction, this specialization is the quantum
dimension. Thus \eqref{eq:1} follows from \cref{thm:2}.

In order to prove \eqref{eq:1}, we may assume $\bw^s$ is a fundamental
weight: there is a torus action on framing vector spaces, and the
induced action on $\fM_{\zeta}(\bv^s,\bw^s)$. Let us suppose
$\bw^s = \Lambda_i+\Lambda_j$ for $i,j\in I^0$ for simplicity. Then
the torus fixed point set is
\begin{equation*}
  \bigsqcup_{\bv^1+\bv^2=\bv^s} \fM_{\zeta}(\bv^1,\Lambda_i)
  \times \fM_{\zeta}(\bv^2,\Lambda_j).
\end{equation*}
As Euler number is equal to the sum of Euler numbers of fixed points
with respect to a torus action, \eqref{eq:1} for $\bw^s$ follows from
\eqref{eq:1} for $\Lambda_i$ and $\Lambda_j$.
This result is compatible with what we mentioned in Introduction:
standard modules are tensor products of $l$-fundamental modules
\cite{VV-std}, as $l$-fundamental modules correspond to the case
$\bw^s$ is a fundamental weight.

A standard module depends also on spectral parameter, which is
specialization homomorphism
$K_{\prod_{i\in I^0} \GL(\bw^s_i)}(\mathrm{pt})\to \CC$. But the
restriction of a standard module to $\mathbf U_q(\fg_{\fin})$ is
independent of the spectral parameter. Hence the spectral parameter is
not relevant for \cref{thm:2}.

\begin{Remark}
  In order to prove \cref{thm:1}, we need to check \eqref{eq:1} only
  when $\bw^s$ is contained in the root lattice. But $\Lambda_i$ does
  not satisfy this condition in general, hence the above reduction
  cannot be performed among $\bw^s$ in the root lattice.
\end{Remark}

\section{Quantum dimensions of standard modules}\label{sec:2}

\subsection{Quantum dimension}

Let $V$ be a finite dimensional representation of $\fg_{\fin}$. The
specialized character
\begin{equation*}
  \left.\ch V\right|_{e^{-\alpha_i}=\exp(\frac{2\pi i}{h^\vee+1})}
\end{equation*}
is called the \emph{quantum dimension} of $V$, denoted by $\dim_q V$,
and was introduced by Andersen~\cite[Def.~3.1]{MR1182414} (see also
Parshall-Wang~\cite{MR1247547}), where
$\zeta = \exp(\frac{2\pi i}{2(h^\vee+1)})$. Note $2\rho$ in
$K_{2\rho} = q^{2\rho}$ in \cite{MR1182414} should be understood as an
element in the dual of the weight lattice by
$\langle 2\rho,\lambda\rangle = (2\rho,\lambda)/(\alpha_0,\alpha_0)$
for a weight $\lambda$. Here $\alpha_0$ is a short root. See
\cite[Lemma~1.1]{MR1247547}.
It was assumed that $\zeta$ is a primitive $\ell$-th root of unity
with odd $\ell$ in \cite{MR1182414}, the definition still makes sense
for our choice. By the Weyl character formula, we have (see
\cite[(3.2)]{MR1182414}, \cite[Th.~1.3]{MR1247547})
\begin{equation*}
  \left.\ch V\right|_{e^{-\alpha_i}=\exp(\frac{2\pi i}{h^\vee+1})}
  = \dim_q V
  = \prod_{\alpha\in\Delta^+_{\fin}} \frac{
    \zeta^{d_\alpha\langle\lambda+\rho,\alpha^\vee\rangle} -
    \zeta^{-d_\alpha\langle\lambda+\rho,\alpha^\vee\rangle}}{
    \zeta^{d_\alpha\langle\rho,\alpha^\vee\rangle} -
    \zeta^{-d_\alpha\langle\rho,\alpha^\vee\rangle}},
\end{equation*}
where $\lambda$ is the highest weight of an irreducible representation
$V = V(\lambda)$, $\Delta^+_{\fin}$ is the set of positive roots of
$\fg_{\fin}$, $d_\alpha\in \{1,2,3\}$ is the square length of $\alpha$
divided by the length of a short root, and $\rho$ is the half sum of
positive roots. (Since we are only considering type $ADE$, we have
$d_\alpha = 1$ for any $\alpha$.)
We have
$d_\alpha\langle\rho,\alpha^\vee\rangle\le \langle\rho,\theta\rangle =
h^\vee$, hence the denominator does not vanish.

Let us introduce $\zeta$-integers by
\begin{equation*}
  [n]_{\zeta} \defeq \frac{\zeta^{n} - \zeta^{-n}}
  {\zeta - \zeta^{-1}},
\end{equation*}
so that
\begin{equation}\label{eq:4}
  \left.\ch V\right|_{e^{-\alpha_i}=\exp(\frac{2\pi i}{h^\vee+1})}
  = \prod_{\alpha\in\Delta^+}
  \frac{[d_\alpha\langle\lambda+\rho,\alpha^\vee\rangle]_{\zeta}}
  {[d_\alpha\langle\rho,\alpha^\vee\rangle]_{\zeta}}.
\end{equation}
We have
\begin{equation*}
  [2(h^\vee+1) + k]_{\zeta} = [k]_{\zeta}
\end{equation*}
as $\zeta^{2(h^\vee+1)} = 1$, as well as
\begin{equation}\label{eq:2}
  [h^\vee+1-k]_{\zeta} = [k]_{\zeta}
\end{equation}
as $\zeta^{h^\vee+1-k} \zeta^{k} = \zeta^{h^\vee+1} = -1$. In
particular, we have $[h^\vee+1]_{\zeta} = 0$. The former is usual
analogy between roots of unity and characteristic $2(h^\vee+1)$. The
latter is a new feature at an even root of unity.

\subsection{type $A$}

Consider type $A_{n-1}$, i.e., $\fg_{\fin} = \algsl(n,\CC)$. We have
$h^\vee = n$.

In type $A_{n-1}$, it is known that the $k$th $l$-fundamental module
is the $k$th fundamental representation of $\SL(n)$.

The quantum dimension of the $k$th fundamental representation
$V(\Lambda_k)$ is
\begin{equation*}
  \dim_q V(\Lambda_k)
  = \begin{bmatrix}
    n \\ k
  \end{bmatrix}_{\zeta}
  = \frac{[n]_{\zeta}[n-1]_{\zeta}\dots [n-k+1]_{\zeta}}
  {[k]_{\zeta} [k-1]_{\zeta} \dots [1]_{\zeta}}.
\end{equation*}
\begin{NB}
In fact, the denominator of \eqref{eq:4} is equal to
\begin{equation*}
  \prod_{1\le i < j\le n} [j-i]_\zeta,
\end{equation*}
as positive roots $\alpha$ are
$\alpha_{i+1} + \dots + \alpha_j$.

On the other hand, from the numerator, we have
\begin{equation*}
  \prod_{1\le i < j\le n} [\lambda_i - \lambda_j + j-i]_\zeta,
\end{equation*}
where $\lambda$ is regarded as a partition $(1^k)$.
Since $\lambda_i - \lambda_j$ is $1$ if $i\le k$, $j>k$ and $0$
otherwise, we have
\begin{equation*}
  \prod_{1\le i<j\le n} \frac{[\lambda_i-\lambda_j+j-i]_\zeta}
  {[j-i]_\zeta} =
  \prod_{i=1}^k \prod_{j=k+1}^n \frac{[j-i+1]_\zeta}{[j-i]_\zeta}
  = \prod_{i=1}^k \frac{[n+1-i]_\zeta}{[k+1-i]_\zeta},
\end{equation*}
as required.
\begin{NB2}
  Older version: the denominator is
  \begin{equation*}
    [n-1]_{\zeta}[n-2]^2_{\zeta}\dots [1]_{\zeta}^{n-1}.
  \end{equation*}
  The numerator is
  \begin{equation*}
    [n]_{\zeta} [n-1]_{\zeta}^2 \dots [n-k+1]_{\zeta}^{k}
    [n-k]_{\zeta}^{k} \dots [k+1]_{\zeta}^{n-k-1} [k]_{\zeta}^{n-k-1}
    \dots [1]_{\zeta}^{n-2},
  \end{equation*}
  where we assume $k\le n-k$. This is checked by considering whether
  $\langle\lambda,\alpha^\vee\rangle$ is $0$ or $1$. In fact, among
  $(n-1)$ coroots $\alpha^\vee$ with
  $\langle\rho,\alpha^\vee\rangle = 1$, we have one $\alpha^\vee$ with
  $\langle\Lambda_k,\alpha^\vee\rangle = 1$. All others have
  $\langle\Lambda_k,\alpha^\vee\rangle = 0$. Thus we get
  $[2]_{\zeta}[1]^{n-2}_{\zeta}$.
  Among $(n-2)$ coroots $\alpha^\vee$ with
  $\langle\rho,\alpha^\vee\rangle = 2$, we have two $\alpha^\vee$ with
  $\langle\Lambda_k,\alpha^\vee\rangle = 1$. All others have
  $\langle\Lambda_k,\alpha^\vee\rangle = 0$. Thus we get
  $[3]_{\zeta}^2 [2]^{n-4}_{\zeta}$. This continues until
  $\langle\rho,\alpha^\vee\rangle = k$. At this last stage we get
  $[k+1]_{\zeta}^k[k]_{\zeta}^{n-2k}$. In total we get
  \begin{equation*}
    [k+1]^k_{\zeta} [k]_{\zeta}^{n-k-1}\cdots
    [2]_{\zeta}^{n-3} [1]^{n-2}_{\zeta}
  \end{equation*}

  Next at $(n-k-1)$ $\alpha^\vee$ with
  $\langle\rho,\alpha^\vee\rangle = k+1$, we get
  $[k+2]_{\zeta}^k [k+1]_{\zeta}^{n-2k-1}$. This will
  continue until $k$ $\alpha^\vee$ with
  $\langle\rho,\alpha^\vee\rangle = n-k$, where we get
  $[n-k]_{\zeta}^{k}$. In this stage we get
  \begin{equation*}
    [n-k+1]_{\zeta} [n-k]_{\zeta}^{k-1} \cdots
    [k+1]_{\zeta}^{n-2k-1}.
  \end{equation*}

  After that, we get
  \begin{equation*}
    [n-k+1]_{\zeta}^{k-1} \dots
    [n-1]^2_{\zeta} [n]_{\zeta},
  \end{equation*}
  as all roots get $1$ by $\langle\Lambda_k,\alpha^\vee\rangle$. In
  total, we get the above.
\end{NB2}%
\end{NB}%
By \eqref{eq:2}, this is equal to $1$, as $[n-i]_\zeta$ cancels with
$[i+1]_\zeta$.

Hence \cref{thm:2} is proved for type $A$.

\begin{NB}
More generally we have
\begin{equation*}
  \dim_q V(\lambda) = \prod_{1\le i<j\le n}
  \frac{[\lambda_i-\lambda_j+j-i]_\zeta}{[j-i]_\zeta}
\end{equation*}
can be expressed by a specialized Schur polynomial.
\end{NB}%

\subsection{type $D$}
Consider type $D_n$, i.e., $\fg_{\fin} = \so(2n,\CC)$. We have
$h^\vee = 2n-2$.

It is known that the $k$th $l$-fundametal module for $1\le k\le n-2$
is isomorphic as an $\so(2n,\CC)$-module to
\begin{equation*}
  \Wedge^k(\CC^{2n}) \oplus \Wedge^{k-2}(\CC^{2n})\oplus \cdots,
\end{equation*}
where $\cdots$ ends as $\Wedge^1(\CC^{2n}) = \CC^{2n}$ if $k$ is odd,
and $\Wedge^0(\CC^{2n}) = \CC$ if $k$ is even. See
\cite[Remark~5.9]{MR1988990}. For $k=n-1$, $n$, the $l$-fundamental
module is isomorphic to the $k$th fundamental representation of
$\so(2n,\CC)$.

As in \cite{MR1153249}, positive roots are
\begin{equation*}
  \{ L_i + L_j \mid i < j\} \sqcup \{ L_i - L_j \mid i < j\},
\end{equation*}
and simple roots are
\begin{equation*}
  \alpha_i = L_i - L_{i+1} \quad (i=1,\dots, n-1),\qquad
  \alpha_n = L_{n-1} + L_{n}
\end{equation*}
with the standard inner product $(L_i,L_j) = \delta_{ij}$.
Fundamental weights are
\begin{equation*}
  \begin{split}
    & \Lambda_i = L_1+ L_2 + \dots + L_i\quad (i=1,\dots, n-2), \\
    & \Lambda_{n-1} = \frac12(L_1 + L_2 + \dots + L_{n-2} + L_{n-1} - L_n),\\
    & \Lambda_n = \frac12(L_1 + L_2 + \dots + L_{n-2} + L_{n-1} + L_n).
  \end{split}
\end{equation*}
Weyl vector $\rho$ is
\begin{equation*}
  \rho = \sum_{i=1}^n (n-i) L_i.
\end{equation*}

We have
\(
   (\Lambda_k, L_1+L_2) = 2
\)
unless $k=1$, $n-1$, $n$. If this holds, we have
\begin{equation*}
  [(\Lambda_k + \rho, L_1+L_2)]_\zeta = [2n-1]_\zeta = 0.
\end{equation*}
Thus
\begin{equation*}
  \dim_q V(\Lambda_k) = 0 \quad \text{if $k\neq 1$, $n-1$, $n$}.
\end{equation*}

We have
\begin{equation*}
  \dim_q V(\Lambda_1) = \frac{[n]_\zeta [2n-2]_\zeta}{[1]_\zeta [n-1]_{\zeta}} = 1,
\end{equation*}
as $[2n-2]_\zeta = [1]_\zeta$, $[n]_\zeta = [n-1]_\zeta$ by \eqref{eq:2}.
 
Next we consider the case $\lambda = \Lambda_{n-1}$. We have
\begin{equation*}
  (\Lambda_{n-1},L_i - L_j) =
  \begin{cases}
    0 & \text{if $j\le n-1$},\\
    1 & \text{if $j=n$},
  \end{cases}
  \qquad
  (\Lambda_{n-1},L_i + L_j) =
  \begin{cases}
    1 & \text{if $j\le n-1$},\\
    0 & \text{if $j=n$}.
  \end{cases}
\end{equation*}
These imply
\begin{equation*}
  \dim_q V(\Lambda_{n-1}) =
  \frac{[2n-2]_\zeta[2n-4]_\zeta \dots [4]_\zeta}
  {[n-1]_\zeta[n-2]_\zeta\dots [3]_\zeta}
  =
  \begin{cases}
    \frac{[2n-2]_\zeta[2n-4]_\zeta \dots [n]_\zeta}
    {[n-1]_\zeta [n-3]_\zeta \dots [3]_\zeta} & \text{if $n$ is even},\\
    \frac{[2n-2]_\zeta[2n-4]_\zeta \dots [n+1]_\zeta}
    {[n-2]_\zeta [n-4]_\zeta \dots [3]_\zeta} & \text{if $n$ is odd}.
  \end{cases}
\end{equation*}
This is $1$ by \eqref{eq:2}. Note $[2n-2]_\zeta = [1]_\zeta = 1$. We
also have $\dim_q V(\Lambda_n) = 1$ as we have a diagram automorphism
$n\leftrightarrow n-1$, $i\leftrightarrow i$ ($i\neq n-1,n$).

In summary,
\begin{equation*}
  \dim_q V(\Lambda_k) =
  \begin{cases}
    1 & \text{if $k=1,n-1, n$},\\
    0 & \text{otherwise}.
  \end{cases}
\end{equation*}
Hence \cref{thm:2} is proved.

\subsection{type $E_6$}
We have $h^\vee = 12$.
The numbering of vertices is
$
\dynkin[edge
  length=.6cm,label]{E}{6}.$ 
  We need to compute quantum dimensions of $V(\Lambda_k)$, as well as
  of $V(\Lambda_1+\Lambda_6)$ since it appears in standard
  representations when $\bw^s$ is the $4$th fundamental
  weight.

Positive roots are
\begin{gather*}
  \alpha_1, \dots, \alpha_6, \alpha_1 + \alpha_3, \alpha_3 + \alpha_4,
  \alpha_4 + \alpha_5, \alpha_5 + \alpha_6, \alpha_2+ \alpha_4,\\
  \alpha_1 + \alpha_3 + \alpha_4, \alpha_2+\alpha_3+\alpha_4,
  \alpha_2+\alpha_4+\alpha_5, \alpha_3+\alpha_4+\alpha_5,
  \alpha_4+\alpha_5+\alpha_6,\\
  \alpha_1+\alpha_2+\alpha_3+\alpha_4,
  \alpha_1+\alpha_3+\alpha_4+\alpha_5,
  \alpha_2+\alpha_3+\alpha_4+\alpha_5,\\
  \alpha_2+\alpha_4+\alpha_5+\alpha_6,
  \alpha_3+\alpha_4+\alpha_5+\alpha_6,\\
  \alpha_1+\alpha_2+\alpha_3+\alpha_4+\alpha_5,
  \alpha_1+\alpha_3+\alpha_4+\alpha_5+\alpha_6,\\
  \alpha_2+\alpha_3+\alpha_4+\alpha_5+\alpha_6,
  \alpha_2+\alpha_3 + 2\alpha_4 + \alpha_5,\\
  \alpha_1+\alpha_2+\alpha_3+\alpha_4+\alpha_5+\alpha_6,
  \alpha_1+\alpha_2+\alpha_3+2\alpha_4+\alpha_5,
  \alpha_2+\alpha_3+2\alpha_4+\alpha_5+\alpha_6,
  \\
  \alpha_1+\alpha_2+2\alpha_3+2\alpha_4+\alpha_5,
  \alpha_1+\alpha_2+\alpha_3+2\alpha_4+\alpha_5+\alpha_6,
  \alpha_2+\alpha_3+2\alpha_4+2\alpha_5+\alpha_6,
  \\
  \alpha_1+\alpha_2+2\alpha_3+2\alpha_4
  +\alpha_5+\alpha_6,
  \alpha_1+\alpha_2+\alpha_3+2\alpha_4
  +2\alpha_5+\alpha_6,
  \\
  \alpha_1+\alpha_2+2\alpha_3+2\alpha_4
  +2\alpha_5+\alpha_6,
  \\
  \alpha_1+\alpha_2+2\alpha_3+3\alpha_4
  +2\alpha_5+\alpha_6,\\
  \alpha_1+2\alpha_2+2\alpha_3+3\alpha_4+2\alpha_5+\alpha_6.
\end{gather*}
Hence the denominator of \eqref{eq:4} is
\begin{equation*}
  [1]_\zeta^6 [2]_\zeta^5 [3]_\zeta^5 [4]_\zeta^5 [5]_\zeta^4
  [6]_\zeta^3 [7]_\zeta^3 [8]_\zeta^2 [9]_\zeta [10]_\zeta [11]_\zeta.
\end{equation*}

We compute the quantum dimension of the first fundamental representation as
\begin{equation*}
  \dim_q V(\Lambda_1) =
  \frac{[2]_\zeta[3]_\zeta[4]_\zeta[5]_\zeta^2[6]_\zeta^2[7]_\zeta^2 [8]_\zeta^2
  [9]_\zeta^2[10]_\zeta[11]_\zeta[12]_\zeta}
{[1]_\zeta[2]_\zeta[3]_\zeta[4]_\zeta^2[5]_\zeta^2[6]_\zeta^2[7]_\zeta^2
  [8]_\zeta^2[9]_\zeta[10]_\zeta [11]_\zeta}
= \frac{[9]_\zeta[12]_\zeta}{[1]_\zeta[4]_\zeta}.
\end{equation*}
\begin{NB}
  the first fundamental representation is $27$-dimensional.
\end{NB}%
This is $1$ by \eqref{eq:2}. We also have $\dim_q V(\Lambda_6)=1$ by
the diagram automorphism.

We have
\begin{equation*}
  \begin{split}
    & (\Lambda_k+\rho,\alpha_1+2\alpha_2+2\alpha_3+3\alpha_4+2\alpha_5+\alpha_6) = 13 \quad\text{for $k=2,3,5$},\\
    & (\Lambda_1+\Lambda_6+\rho,\alpha_1+2\alpha_2+2\alpha_3+3\alpha_4+2\alpha_5+\alpha_6) = 13.
  \end{split}
\end{equation*}
Therefore $\dim_q V(\Lambda_k) = 0$ for $k=2,3,5$,
$\dim_q V(\Lambda_1+\Lambda_6) = 0$.
We also have
\begin{equation*}
  (\Lambda_4+\rho,\alpha_1+\alpha_2+2\alpha_3+3\alpha_4
  +2\alpha_5+\alpha_6) = 13,
\end{equation*}
hence $\dim_q V(\Lambda_4)=0$.
\begin{NB}
    In fact, we have
    \begin{equation*}
        \dim_q V(\Lambda_4) = \frac{
        [5]_\zeta[9]_\zeta^2[10]_\zeta[13]_\zeta[14]_\zeta}
      {[1]_\zeta[2]_\zeta^2[3]_\zeta^2[7]_\zeta}
    \end{equation*}
\end{NB}%

It is known that the $k$th fundamental module is isomorphic as a
$\g_{\fin}$-module to
\begin{equation*}
  \begin{cases}
    V(\Lambda_k) & \text{if $k=1,6$},\\
    V(\Lambda_2)\oplus V(0) & \text{if $k=2$},\\
    V(\Lambda_3)\oplus V(\Lambda_6) & \text{if $k=3$},\\
    V(\Lambda_5)\oplus V(\Lambda_1) & \text{if $k=5$},\\
    V(\Lambda_4)\oplus V(\Lambda_2)^{\oplus 2}
    \oplus V(\Lambda_1+\Lambda_6) \oplus V(0) & \text{if $k=4$}.
  \end{cases}
\end{equation*}
This can be given by using the algorithm \cite{MR2144973}, as we did
for $E_8$ in \cite{Na:E8} with much smaller efforts. Instead, the list
can be found in \cite{MR1436775}, which assumed fermionic formula
conjectured at that time, and proved later in \cite{MR2428305}
\begin{NB}
, based
on many earlier works including \cite{MR1993360}    
\end{NB}%
.

We substitute the above computation of quantum dimensions of various
modules to the above combination.
We find that the combination has the quantum dimension always
$1$. Hence \cref{thm:2} is proved.

\subsection{type $E_7$}
We have $h^\vee = 18$.
The numbering of vertices is
$
\dynkin[edge
  length=.6cm,label]{E}{7}.$ 

We calculate as above, using Sage:
\begin{equation*}
    \dim_q V(\Lambda_7) = \frac{[10]_\zeta[14]_\zeta[18]_\zeta}
    {[1]_\zeta[5]_\zeta[9]_\zeta} = 1
\end{equation*}
and
\begin{equation*}
    \dim_q V(\Lambda_k) = 0\quad \text{if $k\neq 7$}.
\end{equation*}

We have a new pattern:
\begin{equation*}
    \dim_q V(2\Lambda_1) =
    \frac{[12]_\zeta[13]_\zeta[14]_\zeta[15]_\zeta[18]_\zeta[21]_\zeta}{
    [1]_\zeta[2]_\zeta[4]_\zeta[5]_\zeta[6]_\zeta[7]_\zeta}
  = \frac{[21]_\zeta}{[2]_\zeta}.
\end{equation*}
We use \eqref{eq:2} for $k=-2$ this time:
$[21]_\zeta = [-2]_\zeta$. This is $-[2]_\zeta$. Hence
\begin{equation*}
    \dim_q V(2\Lambda_1) = -1.
\end{equation*}
This is the first example of a representation whose quantum dimension
is $-1$. Similarly we have
\begin{equation*}
    \dim_q V(\Lambda_1+\Lambda_7)
    = \frac{[12]_\zeta[14]_\zeta[15]_\zeta[18]_\zeta[20]_\zeta}
    {[1]_\zeta^2 [4]_\zeta [5]_\zeta [7]_\zeta}
    = -1.
\end{equation*}

On the other hand, we get
\begin{equation*}
    \dim_q V(2\Lambda_7) = \dim_q V(\Lambda_1+\Lambda_6) 
    = \dim_q V(\Lambda_2+\Lambda_7) = 0.
\end{equation*}

It is known that the $k$th $l$-fundamental module is isomorphic as a
$\g_{\fin}$-module to
\begin{equation*}
  \begin{cases}
    V(\Lambda_k) & \text{if $k=7$},\\
    V(\Lambda_1)\oplus V(0) & \text{if $k=1$},\\
    V(\Lambda_2)\oplus V(\Lambda_7) & \text{if $k=2$},\\
    V(\Lambda_6)\oplus V(\Lambda_1)\oplus V(0) & \text{if $k=6$},\\
    V(\Lambda_3)\oplus V(\Lambda_6)\oplus V(\Lambda_1)^{\oplus 2}
    \oplus V(0) & \text{if $k=3$},\\
    V(\Lambda_5)\oplus V(\Lambda_2)^{\oplus 2}
    \oplus V(\Lambda_1+\Lambda_7) \oplus V(\Lambda_7)^{\oplus 2}
    & \text{if $k=5$},\\
    \begin{aligned}[t]
    & V(\Lambda_4)\oplus V(\Lambda_3)^{\oplus 3}
    \oplus V(\Lambda_2+\Lambda_7)^{\oplus 2}
    \oplus V(\Lambda_1+\Lambda_6) \oplus V(2\Lambda_1)\\
    & \quad \oplus
    V(\Lambda_6)^{\oplus 4}\oplus V(\Lambda_1)^{\oplus 4}
    \oplus V(2\Lambda_7)\oplus V(0)^{\oplus 2} 
    \end{aligned}
    & \text{if $k=4$}.
  \end{cases}
\end{equation*}
Substituting the above computation, we find that these combinations
have its quantum dimension always $1$: for example, we have
$\dim_q V(2\Lambda_1) + 2\dim_q V(0) = (-1)+2=1$ in the last case.
Hence \cref{thm:2} is proved.

\begin{NB}
  $W(\Lambda_2) = V(\Lambda_2) \oplus q V(\Lambda_7)$,
  
  $W(\Lambda_6) = V(\Lambda_6) \oplus q V(\Lambda_1) \oplus q^2 V(0)$,

  $W(\Lambda_3) = V(\Lambda_2) \oplus q V(\Lambda_6) \oplus (q+q^2)V(\Lambda_1) \oplus q^3 V(0)$

  $W(\Lambda_5) = V(\Lambda_5) \oplus (q+q^2) V(\Lambda_2) \oplus q V(\Lambda_1+\Lambda_7) \oplus (q^2+q^3) V(\Lambda_7)$
\end{NB}%

\subsection{type $E_8$}
We have $h^\vee = 30$.
The numbering of vertices is
$
\dynkin[edge
  length=.6cm,label]{E}{8}.$ 

We have similar patterns:
\begin{equation*}
    \dim_q V(\Lambda_k) = 0 \quad\text{for any $k=1,\dots,8$},
\end{equation*}
\begin{equation*}
    \dim_q V(2\Lambda_8) =
    \frac{[20]_\zeta[21]_\zeta[24]_\zeta[25]_\zeta[30]_\zeta[33]_\zeta}
    {[1]_\zeta[2]_\zeta[6]_\zeta[7]_\zeta[10]_\zeta[11]_\zeta} = -1,
\end{equation*}
\begin{equation*}
    \dim_q V(\Lambda_7+\Lambda_8) =
    \frac{[14]_\zeta[18]_\zeta[20]_\zeta[22]_\zeta[24]_\zeta[25]_\zeta[26]_\zeta
      [30]_\zeta[32]_\zeta[34]_\zeta}
    {[1]_\zeta^2[3]_\zeta[5]_\zeta[6]_\zeta[7]_\zeta[9]_\zeta[11]_\zeta[13]_\zeta
    [17]_\zeta} = 1,
\end{equation*}
\begin{equation*}
    \dim_q V(\Lambda_6+\Lambda_8) =
    \frac{[15]_\zeta[18]_\zeta[20]_\zeta[21]_\zeta[24]_\zeta[25]_\zeta[26]_\zeta
    [27]_\zeta[30]_\zeta[32]_\zeta[33]_\zeta[35]_\zeta}
  {[1]_\zeta^2[2]_\zeta[4]_\zeta^{2}[5]_\zeta[6]_\zeta[7]_\zeta[10]_\zeta
  [11]_\zeta[13]_\zeta[16]_\zeta} = -1,
\end{equation*}
and
\begin{gather*}
    \dim_q V(\Lambda_1+\Lambda_8) = \dim_q V(2\Lambda_1) =
    \dim_q V(\Lambda_2+\Lambda_8) = \dim_q V(\Lambda_1+\Lambda_7) \\
    = \dim_q V(\Lambda_1+2\Lambda_8) = \dim_q V(2\Lambda_7) = \dim_q
    V(3\Lambda_8) = \dim_q V(2\Lambda_1+\Lambda_8)\\
    = \dim_q V(\Lambda_3+\Lambda_8) = \dim_q V(\Lambda_2+\Lambda_7) = 
    \dim_q V(\Lambda_1+\Lambda_6) = \dim_q V(\Lambda_1+\Lambda_2) = 0.
\end{gather*}

The $\g_{\fin}$-module structure of the $k$th $l$-fundamental module
is known. Let us omit \emph{negligible} modules (i.e., those with
$\dim_q = 0$) for brevity. We have
\begin{equation*}
  \begin{cases}
    V(0) & \text{if $k=1,2,7,8$},\\
    V(2\Lambda_8)\oplus V(0)^{\oplus 2} & \text{if $k=3,6$},\\
    V(2\Lambda_8)^{\oplus 5}\oplus V(\Lambda_7+\Lambda_8)^{\oplus 2}
    \oplus V(0)^{\oplus 4} & \text{if $k=5$},\\
    V(\Lambda_6+\Lambda_8)^{\oplus 4}\oplus V(\Lambda_7+\Lambda_8)^{\oplus 18}
    \oplus V(2\Lambda_8)^{\oplus 23} \oplus V(0)^{\oplus 10} 
    & \text{if $k=4$}.
  \end{cases}
\end{equation*}
Rather miraculously we find that all have quantum dimensions $1$. For
example, we calculate the result as $(-4)+18+(-23)+10 = 1$ in the last
case. Hence \cref{thm:2} is proved.

\section{Rationally smoothness}\label{sec:rationally}

As mentioned in Introduction, the original motivation of study of
$\fM_{\zeta^\bullet}(\bv,\bw)$ in \cite{Na-branching} was its relation
between the restriction from $\fg_{\aff}$ to $\fg_{\fin}$.
However the Euler numbers of $\fM_{\zeta^\bullet}(\bv,\bw)$ did not
play any role in \cite{Na-branching}. Intersection cohomology groups of
$\fM_{\zeta^\bullet}(\bv,\bw)$ appeared instead.

Therefore we ask whether $\fM_{\zeta^\bullet}(\bv,\bw)$ is a rational
homology manifold, i.e., its intersection cohomology complex is
quasi-isomorphic to the constant sheaf. See \cite[1.4]{BM} for the
original definition of a rational homology manifold, and its
equivalence to the present one.
This sounds optimistic, but Borho-MacPherson \cite[\S2.3]{BM} showed
that the nilpotent variety is a rational homology manifold.
Also $\fM_0(\bv,\Lambda)$ is a symmetric power of $\CC^2/\Gamma$,
which has only finite quotient singularity, hence is a rational
homology manifold.

\subsection{Example: Hilbert schemes of two points}

Consider
$\fM_{\zeta^\bullet}(2\delta,\Lambda_0) =
\operatorname{Hilb}^2(\CC^2/\Gamma)$. 
As we already mentioned in \cref{ex:Yam}, Yamagishi showed that
$\fM_{\zeta^\bullet}(2\delta,\Lambda_0)$ is locally isomorphic to the
intersection of the nilpotent cone of $\fg_{\fin}$ and the Slodowy
slice to a sub-subregular orbit around $0$-dimensional
strata.
The transversal slice to the bigger stratum is $\CC^2/\Gamma$. It is
also a rational homology manifold.
Combined with the rational smoothness of the nilpotent variety
mentioned above, we obtain the following.

\begin{Corollary}
    $\fM_{\zeta^\bullet}(2\delta,\Lambda_0) = \mathrm{Hilb}^2(\CC^2/\Gamma)$ is a rational homology manifold.
\end{Corollary}

This result can be also proved from \cite{Na-branching}, together with
\cite{Na-qaff} and some Euler number computation. See the argument
below in a simpler situation.

\begin{NB}
  We look at \cref{ex:Yam} for $\bv^s$, $\bw^s$, and compute the
  multiplicities of the trivial representation in standard
  modules.

  Consider type $A_n$ with $n > 2$. We have
  $\bw^s = \Lambda_2 + \Lambda_{n-1}$. The standard module is
  $\Wedge^2 \CC^{n+1}\otimes \Wedge^{n-1} \CC^{n+1}$. It contains the
  trivial representation with multiplicity $1$.

  Consider type $D_n$. For $\bw^s = 2\Lambda_1$, the standard module
  is the tensor square of the vector representation, which is
  $\operatorname{Sym}^2\CC^{2n}\oplus \Wedge^2\CC^{2n}$, and the
  former is the sum of the irreducible representation of the highest
  weight $2\Lambda_1$ and the trivial representation. (See e.g.,
  \cite[Exercise~19.21]{MR1153249}. Therefore the multiplicity of the
  trivial representation is $1$.
\end{NB}%

\subsection{Counter-example}

Consider $A_1^{(1)}$ with $n=4$, i.e.,
$\operatorname{Hilb}^4(\CC^2/(\ZZ/2)) =
\fM_{\zeta^\bullet}(4\delta,\Lambda_0)$. We take
$\bv' = \dynkin[arrows=false,edge length=.6cm,labels={4,2}]{B}{o*}$.
Then $\bv^s = 2$, $\bw^s = 4$. Hence the transversal slice is
$\fM_0(\bv^s,\bw^s)$, which is the nilpotent orbit in
$\mathfrak{sl}(4)$ of type $(2^2)$. It is well known to be a
non-rational homology manifold. For example, we can argue as
follows. By \cite[Th.~15.1.1]{Na-qaff}
$\dim i_0^! \mathrm{IC}(\fM_0(\bv^s,\bw^s))$ gives the multiplicities
of the trivial representation in the standard module for $\bw^s$. The
latter is $4$th tensor power of the natural $2$-dimensional
representation of $\algsl(2)$. The multiplicity of the trivial
representation is $2$.
\begin{NB}
  $(y+y^{-1})^4 = y^4 + 4y^2 + 6 + 4y^{-2} + y^{-4}
  = (y^4 + y^2 + 1 + y^{-2} + y^{-4}) + 3(y^2 + 1 + y^{-2}) + 2$.
\end{NB}%
Hence $\fM_0(\bv^s,\bw^s)$ is \emph{not} a rational homology manifold,
and $\fM_{\zeta^\bullet}(4\delta,\Lambda_0)$ neither.

Alternatively we argue as follows.
We realize this variety as the Coulomb branch of a quiver
gauge gauge theory of type $A_3$ with
$\bv'' = \dynkin[edge length=.6cm,labels={1,2,1}]{A}{3}$,
$\bw'' = \dynkin[edge length=.6cm,labels={0,2,0}]{A}{3}$. It is the
closure of a $\SL(4)[[z]]$-orbit in the affine Grassmannian for
$\SL(4)$, and hence the intersection cohomology is known by geometric
Satake correspondence. We have
\begin{equation*}
  \dim i^!_0\mathrm{IC}(\cM_C(\bv'',\bw'')) = \dim V_0(2\Lambda_2),
\end{equation*}
where $i_0$ is the embedding of the identity element to the affine
Grassmannian. The right hand side is the $0$-weight space in the
representation of $\SL(4)$ with the highest weight $2\Lambda_2$. This
weight space is $2$-dimensional.

In general, if $\fM_0(\bv^s,\bw^s)$ is happen to be the closure of a
$G[[z]]$-orbit in an affine Grassmannian $\Gr_G = G((z))/G[[z]]$
(e.g., it is true in type $A$ by \cite{MR1968260} or the combination
of \cite{2016arXiv160602002N} and \cite{2016arXiv160403625B}), it is a
rational homology manifold if and only if weight spaces of the
corresponding representation are all $1$-dimensional.\footnote{The
  author thanks Dinakar Muthiah for this remark.}

\bibliographystyle{myamsalpha}
\bibliography{nakajima,mybib,coulomb}

\end{document}